\documentclass[12pt,leqno,a4paper,notitlepage]{amsart}
\usepackage{amssymb}

\newtheorem{teo}{Theorem}[section]
\newtheorem{prop}[teo]{Proposition}
\newtheorem{lema}[teo]{Lemma}

\newtheorem{rmk}[teo]{Remark}
\newtheorem{defn}[teo]{Definition}
\newtheorem{stp}[teo]{\newline}
\newenvironment{prova}{\textbf{Proof.}}{\hfill $\quad \Box$}

\newcommand{\rk}{{\mathrm{rk}}}
\newcommand{\Lcal}{\mathcal{L}}
\newcommand{\Ical}{\mathcal{I}}
\newcommand{\Mcal}{\mathcal{M}}
\newcommand{\Ncal}{\mathcal{N}}
\newcommand{\Fcal}{\mathcal{F}}
\newcommand{\Kcal}{\mathcal{K}}
\newcommand{\Dcal}{\mathcal{D}}
\newcommand{\ol}{\overline}
\newcommand{\ox}{\otimes}
\renewcommand{\:}{\colon}
\newcommand{\IP}{{\mathbb P}}
\newcommand{\lra}{\longrightarrow}
\renewcommand{\o}{\mathcal{O}}

\begin{document}

\title[Compactified Jacobians of curves with spine decompositions]{}
\author[Eduardo Esteves]{}

\maketitle

\begin{center}{\large\bf \vskip-2cm Compactified Jacobians of curves\\ 
\vskip0.2cm with spine decompositions}
\end{center}

\

\begin{center} Eduardo Esteves\footnote{
Supported
by CNPq, Processos 301117/04-7 and 470761/06-7,
by CNPq/FAPERJ, Processo E-26/171.174/2003, and by the
Institut Mittag--Leffler (Djursholm, Sweden).}
\end{center}

\

\section{Introduction}

Let $X$ be a curve, that is, a connected, reduced, projective scheme of 
dimension~1 over an algebraically closed field. 
If $X$ is smooth, and nonrational, $X$ can be embedded 
in a canonical Abelian variety, its Jacobian variety $V_X^0$, the moduli scheme 
for invertible sheaves of degree 0 on $X$. If $X$ is singular, one can still 
consider $V_X^0$, but the Jacobian variety is no longer projective and, in general, 
$X$ cannot be embedded in $V_X^0$ nor, of course, in any of its torsors, 
$V_X^d$, the moduli schemes for degree-$d$ invertible sheaves on $X$.  

Compactifications of the $V_X^d$ have been proposed and studied by many 
authors. The first steps were taken by Igusa \cite{Ig} and Mayer and Mumford 
\cite{MM}, but the first compactification was constructed by D'Souza \cite{DS}, 
\cite{DS2}. Later, Altman and Kleiman \cite{AK0}, \cite{AK}, \cite{AK2} gave two 
different constructions of D'Souza's compactification, which work for families of 
integral curves.

If $X$ is reducible, the $V_X^d$ are not even Noetherian. Nevertheless, 
Oda and Seshadri \cite{OS} constructed various 
compactifications of open subschemes of the $V^d_X$, depending on the choice of 
polarizations, when $X$ has at most ordinary nodes for singularities. In our more 
general context, with no conditions imposed on the singularities, 
the first compactifications were constructed by Seshadri 
\cite{sesh}, and the case of families has been treated by Simpson \cite{simpson}. 
It is worth mentioning as well the compactifications by Caporaso \cite{caporaso} 
and Pandharipande \cite{pandha}, constructed directly over the moduli space of 
stable curves. 

In \cite{ed01} there appeared compactifications of open subschemes of the 
$V^d_X$ that had the important property of being fine, that is, of representing a 
functor. Those compactifications, dependent on the choice of a point $P$ 
on the nonsingular locus of $X$, were only shown to be complete. Here we show 
that they are projective, and give sufficient conditions for 
when they are isomorphic 
to Seshadri's compactifications.

More precisely, let $X_1,\dots,X_n$ be the irreducible components of $X$. 
Seshadri's compactifications depend on the choice of a $n$-tuple 
$\mathfrak a=(a_1,\dots,a_n)$ of positive rational numbers, the so-called 
polarization. Given any 
integer $\chi$, Seshadri \cite{sesh} uses Geometric Invariant Theory to construct 
a projective moduli scheme $U_X(\mathfrak a,\chi)$ 
for S-equivalence classes (see \ref{Sequiv}) of 
torsion-free, rank-1 sheaves $I$ (see \ref{setup}) of Euler characteristic 
$\chi$ on $X$ which are $\mathfrak a$-semistable, i.e. such that
$$
\chi(I|_Y/(\text{torsion}))\geq a_Y\chi,\quad\text{where}\quad 
a_Y:=\sum_{X_i\subseteq Y}a_i,
$$
for every subcurve $Y\subseteq X$.

On the other hand, it is shown in \cite{ed01} the existence of a complete 
scheme $J_X^P(\mathfrak a,\chi)$ representing the functor of 
$P$-quasistable sheaves 
on $X$ with respect to $\mathfrak a$, that is, the $\mathfrak a$-semistable 
sheaves $I$ such that the inequality above is strict for every proper subcurve 
$Y\subset X$ containing $P$. Since $J_X^P(\mathfrak a,\chi)$ admits a universal 
sheaf and $U_X(\mathfrak a,\chi)$ corepresents the functor of 
$\mathfrak a$-semistable sheaves, there is a natural map 
$\Psi\: J_X^P(\mathfrak a,\chi)\to U_X(\mathfrak a,\chi)$.

In this paper we show that $J_X^P(\mathfrak a,\chi)$ is projective (see 
Proposition \ref{qstableproj}) and give sufficient conditions for when $\Psi$ is an 
isomorphism: Our Theorem \ref{thm} states that $\Psi$ is an isomorphism if 
$X$ is locally planar and $a_Y\chi$ is an integer only if $Y\ni P$ or $Y$ 
is a spine, i.e. a connected 
subcurve such that $Y\cap\ol{X-Y}$ consists of separating nodes (ordinary nodes 
of $X$ whose removal disconnects it). In particular, this is the case when 
$X$ is of compact type or even treelike.

Most of the statements in the paper can be immediately 
adapted to families of curves, a task 
left to the reader. The need for the work contained in this paper arose 
in \cite{CCE}, where 
Caporaso, Coelho and I construct and study Abel maps.

Briefly, in Section 2 we introduce the schemes $J_X^P(\mathfrak a,\chi)$, 
there called $J^P_E$ (see \ref{seshjac} for the connection), and show they are 
projective. In Section 3, we discuss spines, and show how a curve that decomposes 
in spines is simpler to study. The various technical results obtained in Section 3 
are combined in Section 4 to obtain our main result: sufficient conditions 
for when $\Psi$ is an isomorphism.

Many thanks are due to L. Caporaso and J. Coelho.

\section{Fine compactified Jacobians}

Fix throughout the paper an algebraically closed field $k$, a 
\emph{curve} $X$ over $k$, that is, a connected, reduced, projective 
scheme $X$ of dimension 1 over $k$, and a point $P$ on the nonsingular locus 
of $X$. Also, denote by $X_1,\dots,X_n$ the irreducible components of $X$.

All schemes are assumed locally of finite type over $k$. All points are assumed 
closed, unless stated otherwise.

\begin{stp}\label{setup}\setcounter{equation}{0}\rm
(\emph{Semistable, torsion-free, rank-$1$ sheaves}) 
Let $I$ be a coherent sheaf on $X$. We say that $I$ is: 
\emph{torsion-free} (or \emph{depth-1}) if its associated points
are generic points of $X$; \emph{rank-$1$} if 
$I$ is invertible on a dense open subset of $X$; and \emph{simple} if 
$\text{End}(I)=k$. Each invertible sheaf on $X$ is torsion-free, rank-1 and 
simple. 

A \emph{subcurve} of $X$ is a closed subscheme that is a curve. For each 
subcurve $Y\subseteq X$, let $I_Y$ denote the restriction of 
$I$ to $Y$ modulo torsion or, in other words, the image of the natural map 
$$
I|_Y\lra\bigoplus_{i=1}^m  (I|_Y)_{\xi_i},
$$
where $\xi_1,\dots,\xi_m$ are the generic points of $Y$. If $I$ is a 
torsion-free (resp. rank-1) sheaf on $X$, so is $I_Y$ on $Y$.

Let $E$ be a locally free sheaf on $X$ of rank $r>0$. Set 
$\mu(E):=\deg(E)/r$, the \emph{slope} of $E$. If $\mu(E)\in\text{\bf Z}$, we 
say that $E$ is a \emph{polarization}. For instance, $\o_X$ is a polarization 
of $X$ of slope 0, the \emph{canonical polarization}.

Assume $I$ is torsion-free, rank-1. 
Call $I$ \emph{semistable} (resp. \emph{stable}) 
with respect to $E$ if $\chi(I\ox E)=0$ and
$\chi(I_Y\ox E|_Y)\geq 0$ (resp. $\chi(I_Y\ox E|_Y)>0$) 
for each proper subcurve $Y\subsetneqq X$. Since
\begin{equation}\label{IoxF}
\chi(I_Y\ox E|_Y)=r\chi(I_Y)+\deg(E|_Y)=r(\chi(I_Y)+\mu(E|_Y)),
\end{equation}
where $\mu(E|_Y):=\deg(E|_Y)/r$, the sheaf $I$ is semistable 
with respect to $E$ if and only if 
$\chi(I)=-\mu(E)$, whence $E$ is a polarization, and 
$\chi(I_Y)\geq-\mu(E|_Y)$ for each proper subcurve $Y\subsetneqq X$. 

If $I$ is semistable and $\chi(I_Y\ox E|_Y)>0$ 
for each proper subcurve $Y\subset X$ containing $P$ we say that 
$I$ is \emph{$P$-quasistable}. Of course, if $I$ is stable then $I$ is 
$P$-quasistable.
\end{stp}

\begin{stp}\label{finejac}\setcounter{equation}{0}\rm
(\emph{The fine compactified Jacobians}) There is a scheme $J_X$
parametrizing torsion-free, rank-1,
simple sheaves on the curve $X$; see \cite{ed01} Thm. B, p. 3048.
More precisely, $J_X$ represents 
the functor that associates to each scheme $T$ the set of
$T$-flat coherent sheaves $\Ical$ on $X\times T$ such that
$\Ical|_{X\times t}$ is torsion-free, rank-1 and simple
for each $t\in T$, modulo equivalence $\sim$. We say that such
sheaves $\Ical$ are torsion-free, rank-1 and simple on $X\times T/T$.
We say that two such sheaves $\Ical_1$ and $\Ical_2$
are equivalent, and denote $\Ical_1\sim\Ical_2$,
if there is an invertible sheaf $\Ncal$ on $T$ such that
$\Ical_1\cong\Ical_2\ox p_2^*\Ncal$, where $p_2\:X\times T\to T$ is the
projection map.

If $T$ is a connected scheme, and $\Ical$ is a torsion-free,
rank-1 sheaf on $X\times T/T$, then $\chi=\chi(\Ical|_{X\times t})$
does not depend on the choice of $t\in T$; we say that
$\Ical$ is a sheaf of Euler characteristic $\chi$ on $X\times T/T$. So, 
there is a natural decomposition
$$
J_X=\coprod_{\chi\in\text{\bf Z}}J^\chi_X,
$$
where $J^\chi_X$ is the subscheme of $J_X$
parametrizing sheaves of Euler characteristic $\chi$.

Fix an integer $\chi$. The scheme $J^\chi_X$ is universally closed
over $k$; see \cite{ed01} Thm. 32, (2), p.~3068. However,
in general, $J^\chi_X$ is neither of finite type nor separated.

Let $E$ be a locally free sheaf on $X$ of slope
$-\chi$. By \cite{ed01} Prop. 34, p. 3071, the subschemes
$J^{ss}_E$ (resp. $J^s_E$, resp. $J^P_E$) of $J^\chi_X$
parametrizing simple and semistable (resp. stable, resp. 
$P$-quasistable) sheaves on $X$ with respect to $E$ are open. 
By \cite{ed01} Thm. A, p. 3047,
$J^{ss}_E$ is of finite type and
universally closed, $J^s_E$ is separated and
$J^P_E$ is complete over $k$. We call $J_E^P$ a \emph{fine compactified Jacobian}.
\end{stp}

\begin{stp}\label{ThetaI}\setcounter{equation}{0}\rm
(\emph{Theta divisors}) For each scheme $S$
and $S$-flat coherent sheaf $\Fcal$ on $X\times S$, there
is an associated invertible sheaf $\Dcal(\Fcal)$ on $S$, called
the \emph{determinant of cohomology} of $\Fcal$. If
$\chi(\Fcal|_{X\times s})=0$ for every $s\in S$, then there is an
associated global section $\sigma_{\Fcal}$ of $\Dcal(\Fcal)$, whose
zero scheme parametrizes those $s\in S$ for which
$\Fcal|_{X\times s}$ admits a nonzero global section. For the
construction and basic properties of
$\Dcal(\Fcal)$ and $\sigma_{\Fcal}$; see \cite{ed01} \S 6.1, p. 3076
and Prop. 44, p. 3078. 

Fix an integer $\chi$, and 
recall that $J^\chi_X$ is
the scheme parametrizing simple, torsion-free,
rank-1 sheaves of Euler characteristic $\chi$ on $X$;
see \ref{finejac}. Let
$\Ical$ be a universal sheaf on $X\times J^\chi_X/J^\chi_X$.
To each locally free sheaf 
$E$ on $X$ with $\mu(E)=-\chi$, 
we associate the invertible sheaf $\Lcal_E:=\Dcal(\Ical\ox p_1^*E)$ 
on $J^\chi_X$
and its global section $\theta_E:=\sigma_{\Ical\ox p_1^*E}$,
where $p_1\:X\times J^\chi_X\to X$ is the projection. Recall that
$\Ical$ is unique up to tensoring with an invertible sheaf
from $J^\chi_X$. So, by \cite{ed01} Prop. 44 (3), p. 3078
(the projection property of the determinant
of cohomology), $\Lcal_E$ and $\theta_E$ are well-defined,
modulo isomorphism, and so is the zero scheme $\Theta_E$ of
$\theta_E$. We call $\Theta_E$ a \emph{theta subscheme}. 

Let $\Lcal^n:=\Lcal^{\ox n}|_{J^P_E}$. In \cite{ed01}, we proved that, 
if $n$ is large enough, then $\Lcal^n$ is generated by its global 
sections; and the morphism these sections define,
$$
\Psi^n\: J^P_E \lra \IP(H^0(J^P_E,\Lcal^n)),
$$ 
restricts to an embedding on $J^s_E$. So $J^s_E$ is quasiprojective. 
In general, $\Psi^n$ is not an embedding. Nevertheless, 
$J^P_E$ is projective; see below.
\end{stp}

\begin{prop}\label{qstableproj}
Let $E$ be a polarization on $X$. Then $J^P_E$ is projective.
\end{prop}

\begin{prova} Without loss of generality, assume that $P\in X_1$. Let $r$ 
be the rank of $E$. Let $F$ be any locally free sheaf on $X$ of rank $rn$ 
such that 
\begin{align*}
\deg(F|_{X_1})=&n\deg(E|_{X_1})-(n-1),\\
\deg (F|_{X_i})=&n{\deg (E|_{X_i})}+1\quad\text{for each $i=2,\dots,n$.}
\end{align*}
For instance, $F$ can be constructed from $E^{\oplus n}$ by a sequence 
of elementary transformations centered at one nonsingular point of $X$ on 
each of $X_2,\dots,X_n$ and at $n-1$ nonsingular points of $X$ on $X_1$. 
(These can be chosen to be all equal to $P$.) 

We claim that $J_{E}^P\subseteq J^s_{F}$. Indeed, let $I$ be any 
$P$-quasistable sheaf with respect to $E$. Then $\chi(I\ox E)=0$ and 
$$
r\chi(I_Y)-\deg(E|_Y)\geq 0
$$
for every proper subcurve $Y\subset X$, with 
equality only if $P\not\in Y$.
First, since $E$ and $F$ have the same 
slope, $\chi(I\otimes F)=0$. Second, let $Y$ 
be any proper subcurve of 
$X$, and $m$ its number of irreducible components. If
$P\not\in Y$ then 
$$
\chi(I_Y\otimes F|_Y)=rn\chi(I_Y)+\deg(F|_Y)\geq
-n\deg(E|_Y)+n\deg(E|_Y)+m=m>0.
$$
And if $P\in Y$ then 
$$
\chi(I_Y\otimes F|_Y)=rn\chi(I_Y)+n\deg ({E}|_Y)+m-n\geq n+m-n=m>0.
$$
Either way, $\chi(I_Y\otimes F|_Y)>0$. So, $I$ is stable
with respect to $F$, proving our claim.

Finally, since $J^s_{F}$ is
quasiprojective by \cite{ed01} Thm. C (4), p. 3048,
so is $J_E^P$. But, since
$J_E^P$ is complete, $J_E^P$ is projective.
\end{prova}

\begin{rmk}\rm Keep the notation used in the
proof of Proposition \ref{qstableproj}. We claim that
$J^P_{E}\supseteq J^{ss}_{F}$, and thus $J^P_{E}=J^{ss}_{F}=J^s_{F}$.
Indeed, let $I$ be a semistable sheaf with respect to $F$. 
Since 
$\mu(E)=\mu(F)$, we have $\chi(I\ox E)=0$. Now, let 
$Y$ be a proper subcurve of $X$, 
and $m$ its number of irreducible components. 
If $P\in Y$ then
$$
\chi(I_Y\otimes E|_Y)=r\chi(I_Y)+\deg(E|_Y)=r\chi(I_Y)
+\frac{\deg(F|_Y)-m+n}{n}\geq
\frac{n-m}{n}>0.
$$
On the other hand, suppose $P\not\in Y$. Since $I$ is
semistable with respect to $F$,
$$
\chi(I_Y)\geq -\frac{\deg(F|_Y)}{rn}.
$$
Let $s$ be the smallest nonnegative
integer such that $rn$ divides $\deg(F|_Y)-s$. Then
$$
\chi(I_Y)\geq-\frac{\deg(F|_Y)-s}{rn}.
$$
Now, since $\deg(F|_Y)-m=n\deg(E|_Y)$,
also $s-m$ is divisible by $n$. But
$m\leq n-1$. Then $s\geq m$, and hence
$$
\chi(I_Y\otimes E|_Y)=r\chi(I_Y)+\deg(E|_Y)
\geq-\frac{\deg(F|_Y)-m}{n}+\deg(E|_Y)=0.
$$
Either way, it follows that $I$ is $P$-quasistable with
respect to $E$.
\end{rmk}

\section{Curves with spine decompositions}

\begin{stp}\label{spine}\setcounter{equation}{0}\rm 
(\emph{Spines})
A point $N\in X$ is called a \emph{separating node} if $N$ is an ordinary node 
of $X$ and $X-N$ is not connected. Since $X$ is itself connected, 
$X-N$ would have two connected components. Their closures are called the 
\emph{tails} attached to $N$.

A (connected) subcurve $Y$ of $X$ is called a 
\emph{spine} if every point in $Y\cap\ol{X-Y}$ is a separating node. Then 
each connected component $Z$ of $\ol{X-Y}$ is a tail intersecting 
$Y$ transversally at a single point on the nonsingular loci of $Y$ and $Z$. 

If $Y$ is a union of spines of $X$ then any connected component of 
$Y$ or $\ol{X-Y}$ is a spine. Two spines $Y_1$ and $Y_2$ of $X$ 
with no common component intersect 
transversally at a single point on the smooth loci 
of $Y_1$ and $Y_2$. 

A tuple $\mathfrak Z:=(Z_1,\dots,Z_q)$ of spines $Z_i$ covering $X$ 
with finite pairwise intersection is called a \emph{spine decomposition} 
of $X$. 
\end{stp}

\begin{prop}\label{ctype}\setcounter{equation}{0} Let 
$\mathfrak Z:=(Z_1,\dots,Z_q)$ be a spine decomposition of $X$. 
Then there is an isomorphism
$$
u\:J_X\lra J_{Z_1}\times\cdots\times J_{Z_q}
$$
sending $[I]$ to $([I|_{Z_1}],\dots,[I|_{Z_q}])$. 
Furthermore, for each integer $\chi$, 
$$
u(J_X^\chi)=\bigcup_{\chi_1+\dots+\chi_q=\chi+q-1}
J_{Z_1}^{\chi_1}\times\cdots\times J^{\chi_q}_{Z_q}.
$$
\end{prop}

\begin{prova} The statements are clearly true if $q=1$. Assume $q\geq 2$. 
Since $Z_i$ and $Z_j$ intersect at at most one point, for $i\neq j$, 
there is at least one $Z_i$ which is a tail of $X$, say for $i=1$. 
Set $Y_1:=Z_1$ and $Y_2:=\ol{X-Z_2}$. 
By induction, since $(Z_2,\dots,Z_q)$ is a spine decomposition of $Y_2$, 
it will be enough to show that there is an isomorphism
$$
u'\:J_X\to J_{Y_1}\times J_{Y_2}
$$
sending $[I]$ to the pair $([I|_{Y_1}],[I|_{Y_2}])$ such that 
$$
u'(J_X^\chi)=\bigcup_{\chi_1+\chi_2=\chi+1}
J_{Y_1}^{\chi_1}\times J^{\chi_2}_{Y_2}
$$
for each integer $\chi$.

We need to show that $u'$ is well-defined. First of all, 
a simple torsion-free, rank-1 sheaf $I$ on $X$ 
must be invertible at
$N$, because otherwise $I\cong I_{Y_1}\oplus I_{Y_2}$. So
the restrictions $I|_{Y_1}$ and $I|_{Y_2}$ are 
torsion-free, rank-1 sheaves. The sum of their Euler characteristics is 
$\chi+1$, where $\chi$ is the Euler characteristic of $I$, 
as it follows 
from applying $\chi(\cdot)$ to the natural exact sequence,
$$
0\lra I|_{Y_2}\ox\o_{Y_2}(-N)\lra I\lra I|_{Y_1}\lra 0.
$$

We claim that a torsion-free, rank-1 sheaf $I$ on $X$,
invertible at $N$, is simple if and only if the restrictions
$I|_{Y_1}$ and $I|_{Y_2}$ are simple. Indeed, 
since $N$ is a node, 
and $I$ is invertible at $N$, there is a natural isomorphism
$$
\text{End}(I)\lra \text{End}(I|_{Y_1})\times_{\text{End}(I|_N)}
\text{End}(I|_{Y_2}),
$$
obtained by restriction, and hence 
$$
\dim\text{End}(I)=\dim\text{End}(I|_{Y_1})\cdot
\dim\text{End}(I|_{Y_2}).
$$
So $\dim\text{End}(I)=1$ if and only if 
$\dim\text{End}(I|_{Y_i})=1$ for $i=1$ and $i=2$.

Let $\Ical$ be a universal simple, torsion-free rank-1 sheaf 
on $X\times J_X/J_X$. 
Then, as seen above, $\Ical|_{X\times s}$
is invertible at $N\times s$ for each $s\in J_X$, and hence 
$\Ical$ is invertible along $N\times J_X$. So, for $i=1,2$, the 
restriction $\Ical|_{Y_i\times J_X}$ is also invertible along 
$N\times J_X$. Since $\Ical|_{Y_i\times J_X}$ agrees with
$\Ical$ on $(Y_i-N)\times J_X$, it follows that 
$\Ical|_{Y_i\times J_X}$ is torsion-free, rank-1 on 
$Y_i\times J_X/J_X$. Moreover, since 
$\Ical|_{X\times s}$ is simple, so is
$\Ical|_{Y_i\times s}$ for $i=1,2$, by our claim above, 
for each $s\in J_X$. The pair 
$(\Ical|_{Y_1\times J_X},\Ical|_{Y_2\times J_X})$ defines $u'$.

As for the inverse to $u'$, 
we construct a map
$$
v'\: J_{Y_1}\times J_{Y_2}\lra J_X
$$
as follows: For each $i=1,2$, let $\Ical_i$ 
be a universal simple, torsion-free rank-1 sheaf on 
$Y_i\times J_{Y_i}$. We may assume that 
$\Ical_i|_{N\times J_{Y_i}}$ is trivial (a rigidification). 
Let $\Mcal_i$ be the pullback of $\Ical_i$ to 
$Y_i\times J_{Y_1}\times J_{Y_2}$ under the projection map. 
Since
$$
\Mcal_1|_{N\times J_{Y_1}\times J_{Y_2}}\cong
\o_{N\times J_{Y_1}\times J_{Y_2}}\cong
\Mcal_2|_{N\times J_{Y_1}\times J_{Y_2}},
$$
and since $N$ is a node of $X$, we may glue $\Mcal_1$ and 
$\Mcal_2$ along $N\times J_{Y_1}\times J_{Y_2}$ to obtain a sheaf $\Mcal$ on 
$X\times J_{Y_1}\times J_{Y_2}$, invertible 
along $N\times J_{Y_1}\times J_{Y_2}$, such that 
$\Mcal|_{Y_i\times J_{Y_1}\times J_{Y_2}}\cong\Mcal_i$ for $i=1,2$. Since 
$\Mcal_1$ and $\Mcal_2$ are flat over $J_{Y_1}\times J_{Y_2}$, so is $\Mcal$. 
Also, since 
$\Mcal_i|_{Y_i\times s}$ is simple, torsion-free and rank-1 on $Y_i$ 
for each $i=1,2$, so is $\Mcal|_{X\times s}$, 
for each $s\in J_{Y_1}\times J_{Y_2}$. 
Let $\Mcal$ define $v'$. 

It is not hard to check that $u'$ and $v'$ are indeed inverse to each other.
\end{prova}

\begin{stp}\label{isov}\setcounter{equation}{0}\rm 
(\emph{Spine decomposition}) Let $\mathfrak Z:=(Z_1,\dots,Z_q)$ be a 
spine decomposition of $X$ and 
$\mathfrak c:=(\chi_1,\dots,\chi_q)$ a $q$-tuple of integers. For 
each $\sigma\in S_q$, where $S_q$ is the permutation group of $\{1,\dots,q\}$, 
there is a natural open and closed embedding
$$
v_{\mathfrak Z}^{\mathfrak c,\sigma}\:
J_{Z_1}^{\chi_1}\times \cdots\times J_{Z_q}^{\chi_q}\lra 
J_X^{\chi_1+\cdots+\chi_q}
$$
sending $([I_1],\dots,[I_q])$ to $[I]$, where $I$ is the unique (modulo 
isomorphism) simple, torsion-free, rank-1 sheaf on $X$ such that 
\begin{equation}\label{Iv}
I|_{Z_{\sigma(i)}}\cong I_{\sigma(i)}\ox\o_{Z_{\sigma(i)}}
\Big(\sum_{j=i+1}^q\sum_{N\in Z_{\sigma(i)}\cap Z_{\sigma(j)}}N\Big)
\quad\text{for $i=1,\dots,q$}.
\end{equation}
In fact, $v_{\mathfrak Z}^{\mathfrak c,\sigma}$ is, after translations, 
a restriction of $u^{-1}$ to an open subscheme, where $u$ is the isomorphism 
of Proposition \ref{ctype}. 
To see that $I$ has indeed Euler characteristic 
$\chi_1+\cdots+\chi_q$, just observe that $I$ has a filtration
\begin{equation}\label{filtI}
0=K_0\subsetneqq K_1\subsetneqq\cdots\subsetneqq K_{q-1}\subsetneqq K_q=I
\end{equation}
where $K_j/K_{j-1}\cong I_{\sigma(j)}$ for $j=1,\dots,q$. We set 
$v_{\mathfrak Z}^{\mathfrak c}:=v_{\mathfrak Z}^{\mathfrak c,\text{id}}$.

Let $E$ be a polarization of $X$ such that $\mu(E|_{Z_i})=-\chi_i$ 
for $i=1,\dots,q$. Then 
$$
J_E^{ss}=\bigcup_{\sigma\in S_q} v_{\mathfrak Z}^{\mathfrak c,\sigma}
(J_{E|_{Z_1}}^{ss}\times\cdots\times J_{E|_{Z_q}}^{ss}).
$$
Indeed, a simple torsion-free, rank-1 sheaf $I$ on $X$ satisfying 
\eqref{Iv} is semistable if the successive quotients of 
\eqref{filtI} are semistable. Conversely, let $I$ be a semistable 
sheaf on $X$. By induction, it is enough to prove that $I|_{Z_i}$ is 
semistable for a certain $i$. 
Now, since $\#(Z_i\cap Z_j)\leq 1$ for $i\neq j$, there 
is a spine $Z_j$ that is a tail. If $I|_{Z_j}$ is semistable we 
are done. If not, then $\chi(I|_{Z_j})\geq\chi_j+1$, and hence 
$\chi(I|_{Z'_j})=\chi'_j$, where $Z'_j:=\bigcup_{\ell\neq j}Z_\ell$ and 
$\chi'_j:=\sum_{\ell\neq j}\chi_\ell$, 
implying that $I|_{Z'_j}$ is semistable. Rinse and repeat. 

Notice that, as a consequence, if $q\geq 2$ then $J^s_E=\emptyset$.
\end{stp}

\begin{defn}\label{polint}\setcounter{equation}{0}\rm 
A polarization $E$ of $X$ is called \emph{integer} at a 
(connected) proper subcurve $Y$ of $X$ if $\mu(E|_Y)$ and $\mu(E|_Z)$, 
for every connected component $Z$ of $\ol{X-Y}$, are integers. 
\end{defn}

\begin{prop}\label{Jss=Js}\setcounter{equation}{0} Let 
$E$ be a polarization of $X$. If $E$ is not integer at any 
subcurve of $X$ then $J_E^{ss}=J_E^s$. If $E$ 
is integer only at subcurves of 
$X$ containing $P$ then $J_E^P=J_E^s$.
\end{prop}

\begin{proof} Let $I$ be a semistable sheaf on $X$ which is 
not stable. Then there is a (connected) proper 
subcurve $Y$ of $X$ such that 
$\chi(I_Y)=-\mu(E|_Y)$. If $I$ is $P$-quasistable then $P\not\in Y$. 
We need only show that $\mu(E|_Z)$ is an integer 
for each connected component $Z$ of $\ol{X-Y}$. 

Now, let $K$ denote the kernel of 
the natural surjection $I\to I_Y$. By the additivity of the Euler 
characteristic, $\chi(K\ox E)=0$. Let $Z_1,\dots,Z_q$ be the 
connected components of $\ol{X-Y}$. For each $i=1,\dots,q$, set 
$K_i:=K|_{Z_i}$ and $N_i:=\text{Ker}(K\to K_i)$. Then 
$K=K_1\oplus\cdots\oplus K_q$. 
Moreover, since $N_i$ is a subsheaf of 
$I$, and $I$ is semistable, $\chi(N_i\ox E)\leq 0$, and hence 
$\chi(K_i\ox E)\geq 0$. But
$$
0=\chi(K\ox E)=\chi(K_1\ox E)+\cdots+\chi(K_q\ox E).
$$
So $\chi(K_i\ox E)=0$, and hence 
$\mu(E|_{Z_i})$ is an integer for each $i=1,\dots,q$.
\end{proof}

\begin{stp}\label{Sequiv}\setcounter{equation}{0}\rm
(\emph{The $S$-equivalence}) Let $E$ be a polarization of $X$ and 
$I$ a semistable sheaf on $X$. Then there are (connected) subcurves 
$Z_1,\dots,Z_q$ covering $X$ and a filtration 
$$
0=I_0\subsetneqq I_1\subsetneqq\cdots\subsetneqq I_{q-1}\subsetneqq 
I_q=I
$$
such that the quotient $I_j/I_{j-1}$ is a stable sheaf on $Z_j$ with 
respect to $E|_{Z_j}$ for each $j=1,\dots,q$. The above filtration is 
called a \emph{Jordan--H\"older filtration}. The sheaf $I$ may 
have many Jordan--H\"older filtrations but the collection of 
subcurves $\mathfrak S(I):=\{Z_1,\dots,Z_q\}$ and the 
isomorphism class of the sheaf
$$
\text{Gr}(I):=I_1/I_0\oplus I_2/I_1\oplus\cdots\oplus I_q/I_{q-1}
$$
depend only on $I$, by the Jordan--H\"older theorem. 

Notice that 
also $\text{Gr}(I)$ is semistable, and 
$$
\text{Gr}(I)\cong\bigoplus_{Z\in\mathfrak S(I)}\text{Gr}(I)_Z.
$$ 
We say that two semistable 
sheaves $I$ and $J$ on $X$ are \emph{$S$-equivalent} if 
$\mathfrak S(I)=\mathfrak S(J)$ and $\text{Gr}(I)\cong\text{Gr}(J)$.
\end{stp}

\begin{lema}\label{Gr=}\setcounter{equation}{0} Let $E$ be a 
polarization of $X$, and $I$ and $K$ two 
$S$-equivalent $P$-quasistable sheaves on $X$. If $\mathfrak S(I)$ 
is a collection of spines of $X$, then $I\cong K$.
\end{lema}

\begin{prova} Since the map $u$ in Proposition \ref{ctype} is an isomorphism, 
to show that $I\cong K$, we need only show that 
$I|_Z\cong K|_Z$ for each $Z\in\mathfrak S(I)$.

For each $Z\in\mathcal S(I)$, let $L_Z:=\text{Gr}(I)_Z$, and   
let $\Delta_Z^P$ be the set of points in $Z\cap\ol{X-Z}$ on the 
connected components of $\ol{X-Z}$ not containing $P$. We need only 
show that
\begin{equation}\label{toshowIZ}
I|_Z\cong L_Z\ox\o_Z\Big(\sum_{N\in\Delta^P_Z}N\Big).
\end{equation}
To prove \eqref{toshowIZ}, let
$$
0=I_0\subsetneqq I_1\subsetneqq\cdots\subsetneqq I_{q-1}\subsetneqq 
I_q=I
$$
be a Jordan--H\"older filtration of $I$. For each $j=1,\dots,q$, 
let $Z_j$ be the support of $I_j/I_{j-1}$, and put 
$Y_j:=Z_1\cup\cdots\cup Z_j$. It follows 
from \cite{ed01} Prop. 6, p. 3053, that $P\in Z_1$ and 
that all the $Y_j$ are connected. So, for each $j=1,\dots,q$, 
$$
\Delta^P_{Z_j}=Z_j\cap\ol{X-Y_j}.
$$
On the other hand, since 
$$
I_j\cong I_{j+1}|_{Y_j}\ox\o_{Y_j}\Big(-\sum_{N\in Y_j\cap Z_{j+1}}N\Big)
$$ 
for each $j=1,\dots,q-1$, it follows by induction that
$$
I_j\cong I|_{Y_j}\ox\o_{Y_j}\Big(-\sum_{N\in Y_j\cap\ol{X-Y_j}}N\Big).
$$
for each $j=1,\dots,q$. So
$$
I|_{Z_j}=(I|_{Y_j})|_{Z_j}\cong (I_j)|_{Z_j}\ox
\o_{Z_j}\Big(\sum_{N\in Z_j\cap\ol{X-Y_j}}N\Big)\cong
L_{Z_j}\ox\o_{Z_j}\Big(\sum_{N\in\Delta^P_{Z_j}}N\Big)
$$
for each $j=1,\dots,q$.
\end{prova}

\begin{prop}\label{Gr=cor}\setcounter{equation}{0} 
Let $E$ be a polarization of $X$. Assume $E$ is integer 
only at subcurves of $X$ that are spines or contain $P$. 
Then any two $S$-equivalent 
$P$-quasistable sheaves are isomorphic.
\end{prop}

\begin{proof} For every semistable sheaf $I$ on $X$, the polarization 
$E$ is integer at every $Z\in\mathfrak S(I)$. By hypothesis, every 
$Z\in\mathfrak S(I)$ but the one containing $P$ is a spine. But 
the subcurve $Z\in\mathfrak S(I)$ containing $P$ is also a spine 
because $\ol{X-Z}$ is a union of spines. 
Then apply Lemma~\ref{Gr=}.
\end{proof}

\begin{prop}\label{theta12}\setcounter{equation}{0} 
Let $\mathfrak Z=(Z_1,\dots,Z_q)$ be a spine decomposition of $X$, 
and let 
$\mathfrak c=(\chi_1,\dots,\chi_q)$ be a $q$-tuple of integers. 
Let 
$$
v_{\mathfrak Z}^{\mathfrak c}\: 
J_{Z_1}^{\chi_1}\times \cdots\times J_{Z_q}^{\chi_q}\lra J_X^\chi
$$
be the associated embedding. Then, for 
each locally free sheaf $E$ on $X$ satisfying $\mu(E|_{Z_i})=-\chi_i$ 
for $i=1,\dots,q$, we have 
\begin{equation}\label{theta12eq}
(v_{\mathfrak Z}^{\mathfrak c})^{-1}\Theta_E=
\sum_{i=1}^q J_{Z_1}^{\chi_1}\times\cdots\times 
J_{Z_{i-1}}^{\chi_{i-1}}\times\Theta_{E|_{Z_i}}\times J_{Z_{i+1}}^{\chi_{i+1}}
\times\cdots\times J_{Z_q}^{\chi_q}.
\end{equation}
\end{prop}

{\noindent
(The sum is the subscheme defined by the product of the
sheaves of ideals of the summands.)}

\begin{prova} Set $J:=J_{Z_1}^{\chi_1}\times \cdots\times J_{Z_q}^{\chi_q}$, and 
denote by $q\:X\times J\to X$ the projection. Let $\Ical$ be a universal sheaf 
on $X\times J_X^\chi/J_X^\chi$. Set $v:=v_{\mathfrak Z}^{\mathfrak c}$. 
For each $j=1,\dots,q$, let 
$$
\Ical_j:=(1_X,v)^*\Ical|_{Z_j\times J}\ox\o_{Z_j\times J}
\Big(-\sum_{\ell=j+1}^q\sum_{N\in Z_\ell\cap Z_j}N\times J\Big).
$$
Then $\Ical_j$ is equivalent to $p_j^*\Lcal_j$, where 
$\Lcal_j$ is a universal sheaf on 
$Z_j\times J_{Z_j}^{\chi_j}/J_{Z_j}^{\chi_j}$ and
$$
p_j\:Z_j\times J_{Z_1}^{\chi_1}\times \cdots\times J_{Z_q}^{\chi_q}\lra 
Z_j\times J_{Z_j}^{\chi_j}
$$
is the projection onto the indicated factors. Also, there is a 
filtration
$$
0=\Kcal_0\subsetneqq\Kcal_1\subsetneqq\cdots\subsetneqq \Kcal_{q-1}
\subsetneqq \Kcal_q=(1_X,v)^*\Ical
$$
such that 
$\Kcal_j/\Kcal_{j-1}\cong\Ical_j$ for $j=1,\dots,q$. 

Now, tensor the above filtration with $q^*E$, and take determinants 
of cohomology. From the base-change, functorial, 
projection and 
additive properties of determinants of cohomology and associated
sections, \cite{ed01} Prop.~44, p. 3078, we get
\eqref{theta12eq}.
\end{prova}

\begin{lema}\label{F1F2F}\setcounter{equation}{0}
Let $Z_1$ and $Z_2$ be proper 
complementary
subcurves of $X$. Let $L$ be 
an invertible sheaf 
on $X$. Let $F_1$ and $F_2$ be locally free sheaves on $Z_1$ and $Z_2$,
respectively, of the same rank such that $\det F_i\cong L|_{Z_i}$ for
$i=1,2$. Then there is a locally free sheaf $F$ on
$X$ with $\det F\cong L$ such that $F|_{Z_i}=F_i$ for $i=1,2$.
\end{lema}

\begin{prova} 
Let $r$ be the common rank of $F_1$ and $F_2$.
Fix a very ample invertible sheaf $\o_X(1)$ on $X$.
Let $m$ be an integer such that the twist $F_i(m)$ is generated
by global sections for each $i=1,2$. Then, for each $i=1,2$ there
is an exact sequence of the form
\begin{equation}\label{extFi}
0\to\o_{Z_i}(-m)^{\oplus r-1}\to F_i\to M_i\to 0,
\end{equation}
where $M_i:=(\det F_i)(m(r-1))$. The extension \eqref{extFi}
is represented by an
element $v_i\in H^1(Z_i,M_i(m))^{\oplus r-1}$.
Let $M:=L(m(r-1))$. Since $\det F_i\cong L|_{Z_i}$ for
$i=1,2$, we have $M|_{Z_i}\cong M_i$ for $i=1,2$.
The induced map $M\to M _1\oplus M _2$
is injective with finite length cokernel, and hence
induces a surjection
$$
H^1(X,M(m))^{\oplus r-1}\longrightarrow
H^1(Z_1,M_1(m))^{\oplus r-1}\oplus H^1(Z_2,M_2(m))^{\oplus r-1}.
$$
Let $v\in H^1(X,M(m))^{\oplus r-1}$ mapping to $(v_1,v_2)$. Then
$v$ corresponds to an exact sequence
$$
0\to\o_X(-m)^{\oplus r-1}\to F\to M\to 0
$$
restricting to \eqref{extFi} on $Z_i$ for each $i=1,2$. So 
$F|_{Z_i}=F_i$ for $i=1,2$. Also, $F$ is locally free of rank $r$ with 
$\det F\cong M(-m(r-1))\cong L$.
\end{prova}

\begin{lema}\label{septgt}\setcounter{equation}{0} Let 
$E$ be a polarization of $X$, and $\chi:=-\mu(E)$. Let $I$ be 
a simple, semistable sheaf on $X$ such that 
$\mathfrak S(I)$ is a collection of spines. Then, for 
each $v\in T_{J^\chi_X,[I]}$ nonzero there are an integer $m\geq 2$ and a 
locally free sheaf $F$ on $X$ with $\rk(F)=m\rk(E)$ and 
$\det(F)\cong\det(E)^{\ox m}$ such that $[I]\in\Theta_F$ but 
$v\not\in T_{\Theta_F,[I]}$.
\end{lema}

\begin{prova} Let $\mathfrak Z:=(Z_1,\dots,Z_q)$ be a spine 
decomposition with $\mathfrak S(I)=\{Z_1,\dots,Z_q\}$. For each 
$i=1,\dots,q$, set $E_i:=E|_{Z_i}$ and $I_i:=\text{Gr}(I)_{Z_i}$, and 
put $\chi_i:=\chi(I_i)$ for $i=1,\dots,q$. Notice that 
$\mu(E_i)=-\chi_i$ and that $I_i$ is stable with respect to $E_i$.
Put $\mathfrak c:=(\chi_1,\dots,\chi_q)$. As shown in 
\ref{isov}, up to reordering the elements of $\mathfrak S(I)$, 
we may assume that 
$[I]=v_{\mathfrak Z}^{\mathfrak c}([I_1],\dots,[I_q])$.
Moreover, since $v_{\mathfrak Z}^{\mathfrak c}$ is an open embedding, 
there are $v_i\in T_{J_{Z_i}^{\chi_i},[I_i]}$ for $i=1,\dots,q$ such that 
$d(v_{\mathfrak Z}^{\mathfrak c})_{([I_1],\dots,[I_q])}(v_1+\cdots+v_q)=v$. Since 
$v\neq 0$, there is $i$ such that $v_i\neq 0$. Fix such an $i$.

Since $I_i$ is stable with respect to $E_i$, it follows from 
\cite{edsep} Lemma~12, p. 583, that there are an integer $m\geq 2$ 
and a locally free sheaf $F_i$ on $Z_i$ with 
$\rk(F_i)=m\rk(E)$ and $\det(F_i)\cong\det(E_i)^{\ox m}$ 
such that $[I_i]\in \Theta_{F_i}$ but $v_i\not\in T_{\Theta_{F_i},[I_i]}$. 
On the other hand, for each $j\neq i$, it follows 
from \cite{ed01} Thm.~11, p. 3057, that
there is a locally free sheaf $F_j$ on $Z_j$ with rank $m\rk(E)$ and
determinant $\det(F_j)\cong\det(E_j)^{\ox m}$ such that 
$[I_j]\not\in\Theta_{F_j}$. By a repeated application of 
Lemma \ref{F1F2F}, there is a locally free sheaf $F$ on $X$ with 
$\det(F)\cong\det(E)^{\ox m}$ such that $F|_{Z_i}=F_i$ for $i=1,\dots,q$. 
Now, applying Proposition \ref{theta12}, we get that $F$ satisfies 
the statament of the lemma.
\end{prova}

\begin{prop}\label{spineinfinit}\setcounter{equation}{0} 
Let $E$ be a polarization of $X$, and $\chi:=-\mu(E)$. Assume that 
$E$ is integer only at subcurves of $X$ that are spines 
(resp. are spines or contain $P$). 
Then, for every simple, semistable (resp. $P$-quasistable) 
sheaf $I$ on $X$ and each nonzero $v\in T_{J_X^\chi,[I]}$ there are 
an integer $m\geq 2$ and a locally free sheaf $F$ 
on $X$ with rank $m\rk E$ and determinant $\det(E)^{\ox m}$ 
such that $[I]\in\Theta_F$ but $v\not\in T_{\Theta_F,[I]}$.
\end{prop}

\begin{prova} Just observe that the hypothesis implies 
that $\mathfrak S(I)$ is a collection of spines for every semistable 
sheaf $I$. Then apply Lemma \ref{septgt}.
\end{prova}

\section{Coarse compactified Jacobians}

Throughout this section, let $\chi$ be an integer and 
$\mathfrak a=(a_1,\dots,a_n)$ be a $n$-tuple of rational numbers 
summing up to 1. For each subcurve $Y$ of $X$, set $a_Y:=\sum_{X_i\subseteq Y}a_i$.

\begin{stp}\label{seshjac}\setcounter{equation}{0}\rm
(\emph{The coarse compactified Jacobians}) 
According to Seshadri \cite{sesh} D\'ef. 9 and
Remarques on p. 153, a
torsion-free, rank-1 sheaf
$I$ on $X$ is \emph{$\mathfrak a$-semistable} if
$$
\chi(I_Y)\geq a_Y\chi(I)
$$
for each proper subcurve $Y$ of $X$. In addition, $I$ is called
\emph{$\mathfrak a$-stable} 
if the inequalities are strict. (Seshadri worked 
in higher rank as well, what we will not do here.)

Seshadri's notion of stability is encompassed by ours. In other
words, there is a locally free sheaf $E$ on $X$
such that $\mathfrak a$-semistability (resp. $\mathfrak a$-stability)
for torsion-free, rank-1 sheaves of Euler characteristic $\chi$ 
is equivalent to
semistability (resp. stability) with respect to $E$; see
\cite{edsep} Obs. 13, p. 584. In fact,
any locally free sheaf $E$ on $X$
such that
\begin{equation}\label{Ea}
\mu(E|_{X_i})=-a_i\chi\quad\text{for each $i=1,\dots n$}
\end{equation}
has this property. We let $J_X(\mathfrak a,\chi):=J^{ss}_E$ 
and $J_X^P(\mathfrak a,\chi):=J^P_E$ for any such $E$.

In \cite{sesh} Thm. 15, p. 155, Seshadri constructs 
a projective scheme $U_X(\mathfrak a,\chi)$ 
corepresenting the functor $\text{\bf U}$ 
that associates to each scheme $T$ the set of
torsion-free, rank-1 sheaves $\Ical$ on $X\times T/T$ such that
$\Ical|_{X\times t}$ is $\mathfrak a$-semistable and of Euler 
characteristic $\chi$ for each $t\in T$. The 
points on $U_X(\mathfrak a,\chi)$ are in one-to-one correspondence 
with the $S$-equivalence classes of semistable sheaves. 
For a sketch of this 
construction, see the proof of Proposition~\ref{thetacoarse} below.

Since $J_X(\mathfrak a,\chi)$ represents a functor, 
there exists a universal $\mathfrak a$-semistable 
sheaf of Euler characteristic $\chi$ on 
$X\times J_X(\mathfrak a,\chi)/J_X(\mathfrak a,\chi)$, 
and hence a naturally induced map
\begin{equation}\label{Smap}
\Phi\: J_X(\mathfrak a,\chi)\longrightarrow U_X(\mathfrak a,\chi).
\end{equation}
This map is surjective and its fibers parametrize 
$S$-equivalence classes of simple $\mathfrak a$-semistable sheaves.
\end{stp}

\begin{prop}\label{locplanar}\setcounter{equation}{0}  
If $X$ is locally planar then $U_X(\mathfrak a,\chi)$ is reduced.
\end{prop}

\begin{prova} Fix an ample invertible sheaf 
$\o_X(1)$ on $X$; let $d$ denote its degree. 
Since the family of all $\mathfrak a$-semistable sheaves is bounded, 
there is an integer $t$ such that the twist $I^*(t):=I^*\ox\o_X(t)$ is generated 
by global sections and $h^1(X,I^*(t))=0$ 
for every $\mathfrak a$-semistable sheaf $I$ of Euler characteristic $\chi$ on $X$. 
Then there is an injective homomorphism $\o_X\to I^*(t)$. Taking duals, we 
obtain an injection $I(-t)\to \o_X$, which defines a closed subscheme 
$Y$ of $X$ of length $\ell:=\chi(\o_X)-td-\chi$. Let 
$\text{Hilb}_X^\ell$ be the Hilbert scheme of $X$, parametrizing subschemes of 
length $\ell$. Consider the rational map
$$
\alpha\:\text{Hilb}_X^\ell\dashrightarrow U_X(\mathfrak a,\chi)
$$
sending $[Y]$ to $[\Ical_{Y/X}(t)]$, where $\Ical_{Y/X}$ denotes the 
sheaf of ideals of $Y$ in $X$. The map $\alpha$ is defined on the open locus 
$V$ parametrizing subschemes $Y\subset X$ such that $\Ical_{Y/X}(t)$ is 
$\mathfrak a$-semistable. As shown above $\alpha|_V$ is a surjection. 
By \cite{AIK} Cor. 7, p.~7, since $X$ is 
locally planar, $\text{Hilb}_X^\ell$ is reduced, and hence 
$\alpha(V)\subseteq U_X(\mathfrak a, d)_{\text{red}}$.

Now, let $S$ be a scheme and $\Ical$ a torsion-free, rank-1 sheaf on 
$X\times S/S$ such that $\Ical|_{X\times s}$ is $\mathfrak a$-semistable 
for every $s\in S$. Let $\nu\:S\to U_X(\mathfrak a, d)$ be the induced morphism. 
Since $U_X(\mathfrak a, d)$ corepresents the functor $\text{\bf U}$, it is enough 
to show that $\nu(S)\subseteq U_X(\mathfrak a, d)_{\text{red}}$. This 
is a local condition: we need only show that for each $s\in S$ 
there is an open neighborhood $W\subseteq S$ of $s$ such that 
$\nu(W)\subseteq U_X(\mathfrak a, d)_{\text{red}}$. So, let $s\in S$ and 
set $I:=\Ical|_{X\times s}$. Since $I^*(t)$ is generated 
by global sections, there is an injection $\o_X\to I^*(t)$. Since 
$h^1(X,I^*(t))=0$, up to passing to an open neighborhood of $s$, the injection 
lifts to a homomorphism $\o_{X\times S}\to\Ical^*(t)$, which is 
injective on the fibers of $X\times S/S$. Taking duals, we obtain a homomorphism 
$\Ical^*(-t)\to \o_{X\times S}$, which is again injective on the fibers of 
$X\times S/S$. So, we have a well-defined morphism $\mu\:S\to\text{Hilb}_X^\ell$ 
such that $\mu(S)\subseteq V$ and $\nu=\alpha\circ\mu$. Since $\alpha|_V$ 
factors through $U_X(\mathfrak a, d)_{\text{red}}$, so does $\nu$.
\end{prova}

\begin{prop}\label{thetacoarse}\setcounter{equation}{0} 
Let $E$ be a polarization of $X$ such that $\mu(E|_{X_i})=-a_i\chi$ for 
$i=1,\dots,n$. Let $m$ be an integer greater than $1$ and 
$F$ a locally free sheaf on $X$ 
with $\rk(F)=m\rk(E)$ and $\det(F)\cong\det(E)^{\ox m}$. 
Then there is a 
subscheme $\ol\Theta_F\subseteq U_X(\mathfrak a,\chi)$ 
whose inverse image 
under the natural map 
$\Phi\:J_X(\mathfrak a,\chi)\to U_X(\mathfrak a,\chi)$ is 
the theta subscheme $\Theta_F\subseteq J_X(\mathfrak a,\chi)$.
\end{prop}

\begin{proof} We will need to recall the 
construction of $U_X(\mathfrak a,\chi)$. 
Seshadri fixes an ample invertible sheaf 
$\o_X(1)$ on $X$; let $d$ 
denote its degree. Then he chooses and integer $t$ 
large enough that the twist $I(t):=I\ox\o_X(t)$ 
is generated by
global sections and $h^1(X,I(t))=0$ for every $\mathfrak a$-semistable
torsion-free, rank-1 sheaf $I$ on $X$ of Euler characteristic 
$\chi$. Set $c:= dt+\chi$, and let $S$ 
be the scheme parametrizing $\mathfrak a$-semistable, 
torsion-free, rank-1 quotients
$I$ with Euler characteristic $\chi$ of $\o_X^{\oplus c}(-t)$ 
such that the induced map
$H^0(X,\o_X^{\oplus c})\to H^0(X,I(t))$ is an isomorphism. The
scheme $S$ is an open subscheme of Grothendieck's
scheme of quotients of $\o_X^{\oplus c}(-t)$; so there is a universal
quotient $b\:p_1^*(\o_X^{\oplus c}(-t))\to\Ical$ on $X\times S$,
where $p_1\:X\times S\to X$ is the projection.

Let $G:=\text{SL}(c)$. Then there is an action
$\mu\:G\times S\to S$ given by the natural action of $G$ on
$k^{\oplus c}$, which induces one on $\o_X^{\oplus c}$. 
Notice that, from the description of
$\mu$, and the universal property of the quotient $b$,
there is an isomorphism
\begin{equation}\label{isotau}
(1_X,\mu)^*\Ical\to(1_X,q_2)^*\Ical
\end{equation}
satisfying the cocycle condition, where $q_2\:G\times S\to S$ is
the second projection. (For the statement of the cocycle
condition for invertible sheaves, see
\cite{git} Def. 1.6, p.~30. 
The same statement can be made for all
sheaves.)

Seshadri uses 
Geometric Invariant Theory to show that there is a categorical 
quotient of $S$ under $\tau$, and that the quotient is projective. 
Then it is not difficult to show that this quotient corepresents the 
functor $\text{\bf U}$ described in \ref{seshjac}. 

Fix an isomorphism $\Dcal(E)\to\o_{\text{Spec}(k)}$. Let 
$\Lcal_E(\Ical)$ be the determinant of cohomology of 
$p_1^*E\ox\Ical$ with respect to the second projection 
$p_2\:X\times S \to S$. 
The invertible sheaf $\Lcal_E(\Ical)$ is equipped with a 
$G$-linearization. Indeed, by the functorial properties of the 
determinant of cohomology the isomorphism \eqref{isotau} 
induces an isomorphism
$\mu^*\Lcal_E(\Ical)\to q_2^*\Lcal_E(\Ical)$ satisfying the cocycle 
condition. Of course, $\Lcal_E(\Ical)^{\ox\ell}$ comes 
with an induced $G$-linearization for each integer $\ell$.

For each positive 
integer $\ell$, each locally free sheaf $H$ on $X$ with rank 
$\ell\rk(E)$, and each isomorphisms
$\det(H)\to\det(E)^{\ox\ell}$ and $\Dcal(H)\to\o_{\text{Spec}(k)}$, we 
produce a global section $\theta_H(E,\Ical)$ of 
$\Lcal_E(\Ical)^{\ox\ell}$ as follows: First of all, 
using the identification $\det(E^{\oplus\ell})=\det(E)^{\ox\ell}$, 
the isomorphism $\det(H)\to\det(E)^{\ox\ell}$ induces an isomorphism,
\begin{equation}\label{DIFG}
\tau\:\Dcal(\Ical\ox p_1^*H)\ox\Dcal(E)^{\ox\ell}
\longrightarrow\Dcal(\Ical\ox p_1^*(E^{\oplus\ell}))\ox\Dcal(H);
\end{equation}
see \cite{ed01}, Lemma 46, p. 3082. Combining $\tau$ with the 
isomorphisms $\Dcal(E)\to\o_{\text{Spec}(k)}$ and 
$\Dcal(H)\to\o_{\text{Spec}(k)}$, and the identification 
$\Dcal(\Ical\ox p_1^*(E^{\oplus\ell}))=\Dcal(\Ical\ox p_1^*E)^{\ox\ell}$, 
we get an isomorphism 
$\Dcal(\Ical\ox p_1^*H)\to\Lcal_E(\Ical)^{\ox\ell}$. Then we let 
$\theta_H(E,\Ical)$ be the section corresponding to 
$\sigma_{\Ical\ox p_1^*H}$ under this isomorphism. By the 
functorial properties of $\mu$, the section $\theta_H(E,\Ical)$ of 
$\Lcal_E(\Ical)^{\ox\ell}$ is $G$-invariant, that is, 
$\mu^*\theta_H(E,\Ical)$ is carried to $p_2^*\theta_H(E,\Ical)$ 
under the $G$-linearization of $\Lcal_E(\Ical)^{\ox\ell}$ mentioned 
above.

Set $F_0:=F$.
It follows from \cite{ed01} Thm. 11, p. 3057,
that there are locally free sheaves $F_1,\dots,F_r$ on $X$ of rank
$m\text{rk}(E)$ and
determinant $(\det E)^{\ox m}$ such that, for each 
$\mathfrak a$-semistable
torsion-free, rank-1 sheaf $I$ on $X$ of Euler characteristic 
$\chi$, there is
$j$ such that $h^0(X,I\ox F_j)=0$.
Fix isomorphisms
$\det F_j\to(\det E)^{\ox m}$ and $\Dcal(F_j)\to\o_{\text{Spec}(k)}$ for
each $j=0,\dots,r$, and consider the well-defined $G$-invariant
map
$$
(\theta_{F_0}(E,\Ical),\dots,\theta_{F_n}(E,\Ical))\:S\longrightarrow
\text{\bf P}^r.
$$
Since $U_X(a,d)$ is the categorical quotient of $S$ by $G$, the
above map factors through the quotient map $\pi\:S\to U_X(a,d)$.
It follows that there is a subscheme $\ol{\Theta}_F$ of
$U_X(\mathfrak a,d)$ such that 
$\pi^{-1}(\ol{\Theta}_F)$ is the zero scheme
$Z_F$ of $\theta_F(E,\Ical)$, or of $\sigma_{\Ical\ox p_1^*F}$.

We claim
that $\Phi^{-1}(\ol{\Theta}_F)=\Theta_F$. Indeed,
choose a universal sheaf on $X\times J^{ss}_E/J^{ss}_E$; call it $\Ncal$. Let
$z\in J_E^{ss}$ and $I:=\Ncal|_{X\times z}$. Since
$I$ is $\mathfrak a$-semistable of Euler characteristic $\chi$,
the twist $I(t)$ is globally spanned, $h^1(X,I(t))=0$ and 
$h^0(X,I(t))=c$, and there
is a surjection $\o_X^{\oplus c}(-t)\to I$ inducing an isomorphism 
$H^0(X,\o_X^{\oplus c})\to H^0(X,I(t))$. The same
applies to all points $z'$ on a neighborhood
$W\subseteq J^{ss}_E$ of
$z$. So there is a map
$\lambda\:W\to S$ such that
$(1_X,\lambda)^*\Ical\cong\Ncal|_{X\times W}$. Using the base-change
property of the determinant of cohomology and its associated global
section,
it follows that
$\lambda^{-1}(Z_F)=\Theta_F\cap W$. Since
$\Phi=\pi\circ\lambda$, and $z$ was any point of $J^{ss}_E$,
it follows that
$\Phi^{-1}(\ol{\Theta}_F)=\Theta_F$,
as claimed.
\end{proof}

\begin{teo}\label{thm}\setcounter{equation}{0} 
Assume that every subcurve $Y$ of $X$ with $a_Y\chi\in\text{\bf Z}$ 
is a spine or 
contains $P$. Then the natural map 
$\Phi\:J_X(\mathfrak a,\chi)\to U_X(\mathfrak a,\chi)$ 
restricts to a bijective closed embedding
$\Psi\:J_X^P(\mathfrak a,\chi)\to U_X(\mathfrak a,\chi)$. Furthermore, if 
$X$ is locally planar then $\Psi$ is an isomorphism.
\end{teo}

\begin{prova} By Proposition \ref{Gr=cor}, the restriction 
$\Psi:=\Phi|_{J_X^P(\mathfrak a,\chi)}$ is injective. It is also 
surjective, by \cite{ed01} Thm. 7, p. 3054. In addition, 
by Propositions \ref{spineinfinit} and \ref{thetacoarse}, it 
is an immersion. Since $J_X^P(\mathfrak a,\chi)$ is complete, 
$\Psi$ is proper, and thus $\Psi$ is a bijective embedding. If 
$X$ is locally planar then $U_X(\mathfrak a,\chi)$ is reduced by 
Proposition \ref{locplanar}, and hence $\Psi$ is an isomorphism. 
\end{prova}

\end{document}